\documentclass[11pt,a4paper]{article}
\usepackage{graphicx}
\usepackage{color}
\usepackage{amsmath,amsthm,amsfonts}

\newcommand\R{\mathbb{R}}

\newcommand\Z{\mathbb{Z}}

\newcommand\calH{\mathcal{H}}

\newcommand\calZ{\mathcal{Z}}

\newcommand\weakto{\mathop{\rightharpoonup}}
\newcommand\supp{\mathop{\mathrm{supp}}}

\newcommand\dist{{\mathrm{dist}}}

\newcommand\Id{{\mathrm{Id}}}

\newcommand\loc{{\mathrm{loc}}}
\newcommand\eps{{\varepsilon}}

\newcommand{\LM}[1]{\hbox{\vrule width.2pt \vbox to#1pt{\vfill \hrule width#1pt height.2pt}}}
\newcommand{\LL}{{\mathchoice
{\,\LM7\,}{\,\LM7\,}{\,\LM5\,}{\,\LM{3.35}\,}}}
\newtheorem{theorem}{Theorem}[section]
\newtheorem*{theorem*}{Theorem}
\newtheorem{lemma}[theorem]{Lemma}

\newtheorem{remark}[theorem]{Remark}
\newtheorem{corollary}[theorem]{Corollary}

\numberwithin{equation}{section}
\def\normal{\hat\nu}
\def\tildenormal{\hat\nu_*}
\def\e{\varepsilon}
\def\xzero{y}
\usepackage{fancyhdr}
\begin{document}
\begin{center}
{ \LARGE
Density of polyhedral partitions\\[5mm]}
Andrea Braides$^1$, Sergio Conti$^{2}$, Adriana Garroni$^{3}$
\\[2mm]
{\em $^1$ Dipartimento di Matematica, Universit\`a di Roma Tor Vergata,\\
00133 Roma, Italy}\\[1mm]
{\em $^2$ Institut f\"ur Angewandte Mathematik,
Universit\"at Bonn\\ 53115 Bonn, Germany }\\
{\em $^3$ Dipartimento di Matematica, Sapienza, Universit\`a di Roma\\
00185 Roma, Italy}\\[1mm]
\bigskip

\begin{minipage}[c]{0.8\textwidth}
\centerline{\bf Abstract}
We prove the density of polyhedral partitions in the set of finite Caccioppoli partitions. Precisely, we consider
a decomposition $u$ of a bounded Lipschitz set $\Omega\subset\R^n$ into finitely many subsets of finite perimeter, which can be identified with a function in $SBV_{\rm loc}(\Omega;{\cal Z})$ with ${\cal Z}\subset \R^N$ a finite set of parameters. 
For all $\e>0$ we prove that such a $u$ is $\eps$-close to a small deformation of a  polyhedral decomposition $v_\eps$,
in the sense that there is a $C^1$ diffeomorphism $f_\eps:\R^n\to\R^n$ which is $\eps$-close to the identity
and such that $u\circ f_\eps-v_\eps$ is $\eps$-small in the strong $BV$ norm.
This implies that the energy of $u$ is close to that of $v_\eps$ for a large class of energies defined on partitions. Such type of approximations are 
very useful in order to simplify computations in the estimates of $\Gamma$-limits. 
\end{minipage}
\end{center}

\section{Introduction}
Besides their theoretical interest, approximation results have a great technical importance in 
the treatment of variational problems; in particular, in the computation of 
$\Gamma$-limits for varying energies.
The density of piecewise-affine maps in Sobolev spaces, 
for example, often allows computations for integral energies to be performed only in the simplified 
setting of maps with constant gradient. Similarly, the approximation of sets of finite
perimeter by polyhedral sets, which sometimes is taken as the definition of sets of finite perimeter
itself, allows to reduce problems involving surface energies to the case of a planar interface. 
The use of approximation theorems for the computation of $\Gamma$-limits
is not strictly necessary, since more abstract
integral-representation theorems can be used, whose application though is often 
quite technical. The computation is simpler if representation formulas are available such as relaxation 
or homogenization formulas. Indeed, in that case it is easier to prove a lower bound for 
a $\Gamma$-limit by the blow-up technique elaborated by Fonseca and M\"uller \cite{FonsecaMueller1993}. Approximation results are crucial to reduce the proof of the upper bound
to simpler functions for which recovery sequences are suggested by the representation
formulas themselves.

In multi-phase problems, i.e.,  for interfacial problems when more than two sets
are involved, the proper variational setting is that of partitions into sets of finite
perimeter, or Caccioppoli partitions, for which a theory of relaxation and 
$\Gamma$-convergence has been first developed by Ambrosio and Braides \cite{AmbrosioBraides1990a,AmbrosioBraides1990b}.
The study of Caccioppoli partitions is also a fundamental step in the analysis of 
free-discontinuity problems defined on (special) functions of bounded variation,
since lower-semicontinuity conditions and representation formulas for the latter can 
be often deduced from those for partitions. In that spirit, integral-representation theorems 
for partitions have been proved by Braides and Chiad\`o Piat \cite{BraidesChiadopiat1996} and Bouchitt\'e 
et al.~\cite{BouchitteFonsecaLeoniMascarenas2002}.

The scope of this paper is to fill a gap that seemingly exists in the 
treatment of problems on Caccioppoli partitions, namely the existence of
approximations by polyhedral sets. This is a widely expected result, so
much expected that sometimes it is mistakenly referred to as proved in some
reference text. Conversely, its non-availability makes it more complicated 
to obtain homogenization results even when formulas are available (see 
e.g.~the recent work by Braides and Cicalese \cite{BraidesCicalese2015}).
A  ``dual'' result for systems of rectifiable curves has recently been proved by 
Conti et al.~\cite{ContiGarroniMassaccesi2015} and used to show 
convergence of linear elasticity to a dislocation model \cite{ContiGarroniOrtiz2015}.
In a two-dimensional setting one can use that approximation result to obtain 
polyhedral partitions by considering boundaries of sets
as rectifiable curves, and as such it has been recently used to the study
of systems of chiral molecules \cite{BraidesGarroniPalombaro2015}.

It must be noted that the method usually followed
to obtain approximating sets for a single Caccioppoli set 
cannot be used for partitions. Indeed, for
a single set of finite perimeter $E\subset\R^n$ we can use that the characteristic
function $u:=\chi_E$ is by definition a function with bounded variation; hence, by a
mollification argument it can be approximated by smooth functions $u_\rho$
and this mollification process does not increase the corresponding variation.
Approximating sets are then obtained by taking  super-level sets
of the form $E_\rho:=\{x: u_\rho(x)>c_\rho\}$.
By Sard's theorem the set $E_\rho$ is smooth for almost all values of 
$c_\rho$, and by the coarea formula 
 $c_\rho$  can be chosen so that the boundary of $E_\rho$ is not larger than
 the boundary of $E$. Finally, polyhedral sets are obtained by triangulation
using the
smoothness of $E_\rho$. Such a simple argument cannot be repeated if
we have a partition. Indeed, identify such a partition $(E_1,\ldots, E_N)$
in ${\mathbb R}^n$
with a $BV$-function by setting $u:=\sum_j a_j \chi_{E_j}$, for suitable 
labelling parameters $a_j$. If we choose $a_j$ real numbers, the
process outlined above will require the choice of more superlevel sets $\{u_\rho>c_\rho^j\}$,
which will introduce artificial interfaces. To picture this situation, think of
having a partition into three sets of finite perimeter and choose as labels
the numbers $a_j:=j$. Then in the process above we will have two approximating
sets $E_\rho^1:=\{x: c_\rho^1< u_\rho(x)\le c_\rho^2\}$ and 
$E_\rho^2:=\{x:  u_\rho(x)> c_\rho^2\}$ with $1<c_\rho^1< 2< c_\rho^2<3$
and the interface between the set $E_1$ and $E_3$ will be approximated 
by a double interface: one between $E^1_\rho$ and $E^2_\rho$ and another one between $E^2_\rho$ and $E^3_\rho$. Although
these approximations weakly converge to the original partition, the total length
of the surface has doubled and the energy of the partitions will not converge.
If otherwise we label the sets with $a_j$ in some higher-dimensional
${\mathbb R}^m$ then the use of the coarea formula is not possible.
It is then necessary, as is done in the proof of integral-representation results, 
to make a finer use of the structure of boundaries of sets of finite perimeter.

In our construction we use the fact that essential boundaries between sets of finite
perimeter are contained in $C^1$ hypersurfaces that can be locally deformed 
onto portions of hyperplanes.
By a covering argument we can thus transform most of the interfaces with
a small deformation into open subsets of a finite system of hyperplanes,
which can in turn  be approximated by polyhedral sets. We finally introduce
a decomposition of the ambient space into a system of small polyhedra whose
boundaries contain the above-mentioned lower-dimensional polyhedral sets, and define
a Caccioppoli partition by choosing the majority phase (i.e., the label corresponding to
the set with the largest measure) on each of the small polyhedra. 
This finally gives the desired approximating sets. 

A scalar version of this result is proven in \cite[Th. 4.2.20]{Federer1969} and then refined in 
\cite{AmarDecicco2005} and \cite{Quentindegromard2008}.
The vector-valued approximation, however, does not follow from the scalar one working componentwise since
approximation of the energy requires to choose a single deformation $f$ for all components.
Approximation of vector-valued  $SBV^p$ functions was studied in \cite{CortesaniToader1999,KristensenRindler2016}, but
the case of partitions does not seem to follow directly from the arguments therein, which introduce
large gradients in small regions. The vectorial case of $SBD^p$ functions was addressed for $p=2$ in 
 \cite{Chambolle2004,Chambolle2005,Iurlano2014}.

\section{The density result and its proof}
We will consider partitions of an open set $\Omega\subset\R^n$ into $N$ sets of finite perimeter
$(E_1,\ldots, E_N)$ for some $N\ge 1$.
We say that a set $\Sigma\subset\Omega$ is {\em polyhedral} if there is a finite
number of $n-1$-dimensional simplexes $T_1,\dots, T_M\subset \R^n$ such that
$\Sigma$ coincides, up to $\calH^{n-1}$-null sets, with $\bigcup_{j=1}^M T_j\cap \Omega$.
We are interested in showing that for a general partition $(E_1,\ldots, E_N)$
there exist polyhedral approximations; i.e., partitions $(E^j_1,\ldots,E^j_N)$ of $\Omega$ into sets
whose boundaries are polyhedral, such that for all $k\in\{1,\ldots,N\}$ we have
$|E^j_k\triangle E_k|\to 0$ and $\calH^{n-1}(\partial E^j_k)\to \calH^{n-1}(\partial E_k)$
as $j\to+\infty$ (where in the last formula $\partial E_k$ denotes the reduced boundary of $E_k$),
and the normal to $\partial E_k^j$ converges in a suitable sense to the normal to $\partial E_k$ (see Corollary~\ref{Corollario}).

It will be handy to use a finite set $\calZ:=\{z_1,\ldots, z_N\}\subset \R^N$ as a set of labels 
for the different phases, and identify each partition with the function 
$u:\Omega\to\R^N$ given by $u(x)=z_k$ on $E_k$. In this way the 
set of partitions into $N$ sets of finite perimeter is identified with a subset of the space $SBV_\loc(\Omega;\calZ)$ 
(see \cite{AmbrosioBraides1990a,AmbrosioBraides1990b}). Note that 
our results will be independent of the labelling, but the latter allows to
state the convergence of boundaries of sets as a convergence of the
derivatives of functions.

We recall that a function $u\in SBV_\loc(\Omega;\calZ)$, for $\Omega\subset\R^n$ open, 
has the property that its distributional derivative is a bounded measure of the form
$Du=[u]\otimes\nu\calH^{n-1}\LL J_u$. Here $J_u\subset\Omega$ is a $n-1$-rectifiable set,
called the jump set of $u$, the unit vector $\nu:J_u\to S^{n-1}$ is the normal to $J_u$, and $[u]:=(u^+-u^-)$ is the jump of $u$, where $u^+$ and $u^-:J_u\to\calZ$ are the traces of $u$ on the two
sides of $J_u$, which are $\calH^{n-1}\LL J_u$-measurable. 
We use the notation $\mu\LL E$ for 
the restriction of a measure $\mu$ to a $\mu$-measurable set $E$, defined by $(\mu\LL E)(A):=\mu(E\cap A)$.

\bigskip
The main result of this paper is the following approximation statement.
\begin{theorem}\label{theorem} Let
  $\calZ\subset\R^N$ be finite, let $u\in SBV_\loc(\Omega;\calZ)$ with $|Du|(\Omega)<\infty$, and let $\Omega\subset\R^n$ be a Lipschitz set with $\partial\Omega$ compact. Then there is a sequence $u_j\in SBV_\loc(\Omega;\calZ)$ such that
 $J_{u_j}$ is polyhedral,
 $u_j\to u$ in $L^1_\loc(\Omega;\calZ)$ and $Du_j\weakto Du$ as measures,  and there are bijective maps $f_j\in C^1(\R^n;\R^n)$, with inverse also in $C^1$,
 which converge strongly in $W^{1,\infty}(\R^n;\R^n)$ to the identity map such that
 $|D (u\circ f_j)-Du_j|(\Omega)\to0$.
 \end{theorem}

We remark that $u\circ f_j$ is defined on the set $f_j^{-1}(\Omega)$, and so is the measure $D (u\circ f_j)$, which is then
implicitly extended by zero to the rest of $\R^n$.
In particular, 
$$
|D (u\circ f_j)-Du_j|(\Omega)=|D (u\circ f_j)-Du_j|(\Omega\cap f_j^{-1}(\Omega))+
|Du_j|(\Omega\setminus f_j^{-1}(\Omega)).
$$

The rest of this paper contains the proof of Theorem \ref{theorem}. We shall first (Theorem \ref{theoremRn}) prove the analogous statement for functions defined on $\R^n$, 
and then (Lemma \ref{lemmaext}) give an extension argument to deal with general domains.
We remark that the assumption that $\partial\Omega$ is compact is only used in constructing the extension, so that our result can be extended immediately to
some other unbounded sets, such as, for example, the half space.

\begin{theorem}\label{theoremRn}
 Let
  $\calZ\subset\R^N$ be finite,  and let $u\in SBV_\loc(\R^n;\calZ)$ with $|Du|(\R^n)<\infty$. Then there is a sequence $u_j\in SBV_\loc(\R^n;\calZ)$ such that
 $J_{u_j}$ is polyhedral,
 $u_j\to u$ in $L^1_\loc(\R^n;\calZ)$, $Du_j\weakto Du$ as measures,  and there are bijective maps $f_j\in C^1(\R^n;\R^n)$, with inverse also in $C^1$,
 which converge strongly in $W^{1,\infty}(\R^n;\R^n)$ to the identity map such that
 $|D (u\circ f_j)-Du_j|(\R^n)\to0$.
\end{theorem}
The proof of Theorem \ref{theoremRn}
relies on a deformation
argument allowed by the rectifiability of $J_u$. We recall that the latter means that
$J_u$ coincides, up to an $\calH^{n-1}$-null set, with a Borel subset of the union of countably many $C^1$ surfaces
 \cite[Sect. 2.9]{AmbrosioFP}.
Furthermore, in this characterization one also has that, for $\calH^{n-1}$-almost all points $y\in J_u$, denoting by $M_y$ the $C^1$ surface containing $y$,
the vector
$\nu(y)$ is the normal in $y$ to the  surface $M_y$
and 
\begin{equation}\label{blow-up1}
 \lim_{\rho\to0}
\frac{1}{\rho^{n-1}}  \calH^{n-1}( (J_u \triangle M_y) \cap B_\rho(y)) =0,
 \end{equation}
where $B_\rho(y)$ is the open ball of radius $\rho$ centered in $y$.
The measurability of  the traces $u^\pm(y)$ and the finiteness of $\calZ$ imply, via the Lebesgue point theorem,
that the traces are locally approximately constant, in the sense that
\begin{equation}\label{blow-up2}
\lim_{\rho\to0}
\frac{1}{\rho^{n-1}}
\calH^{n-1}(\{x\in J_u\cap B_\rho(y):u^+(x)\ne u^+(y)\})=0
\end{equation}
for $\calH^{n-1}$-almost every $y\in J_u$.
We refer to \cite{AmbrosioFP} for a more detailed treatment of these concepts.
The idea of the proof is to cover most of the jump set of $u$ by disjoint balls, such that in each of them
the jump set is an (almost flat) $C^1$ graph (see Step 2). In each of the balls the jump set can then be explicitly deformed
into a plane, up to an interpolation region (see Step 1). 
 \begin{proof}
\emph{Step $1$. We perform a local construction around $\calH^{n-1}$-almost all points of the jump set.}

 Fix $\eps\in(0,1)$, whose value will be chosen below.
 Assume that $\xzero\in J_u$ has the following properties:
 there are $g=g_{\xzero}\in C^1(\R^{n-1})$, $r=r_y>0$ and an affine isometry $I_{\xzero}:\R^n\to\R^n$,
 $I_{\xzero}(x)=Q_{\xzero}x+b_{\xzero}$,
satisfying $g(0)=0$, $Dg(0)=0$,
 \begin{equation}\label{eqpropg1}
  \calH^{n-1}((I_{\xzero}J_u) \triangle \{(x', g(x')): x'\in B'_r\}) <  \eps r^{n-1},
 \end{equation}
where $B'_r$ denotes the $n-1$-dimensional ball of radius $r$ centered in $0$,
\begin{equation}\label{eqpropg2}
\calH^{n-1}(\{J_u\cap B_r(\xzero):u^+(x)\ne u^+(\xzero)\})<\eps r^{n-1}
\end{equation}
and the same for $u^-$. Since we chose the isometry $I_{\xzero}$ to make $Dg(0)=0$,
choosing $r$ sufficiently small we can ensure that  additionally
$|Dg| \le\eps^2$ in $B_r'$, which in turn implies $|g|\le \eps^2r$ in $B_r'$.
By  (\ref{blow-up1}) and (\ref{blow-up2})
$\calH^{n-1}$-almost every $\xzero\in J_u$ 
has the  properties above.


We fix $\psi\in C^1_c(B_r(\xzero);[0,1])$ such that $\psi=1$ on $B_{(1-\eps) r}(\xzero)$ 
and $\|D\psi\|_\infty \le 2/(\eps r)$ and
define $f:\R^n\to\R^n$ by 
\begin{equation*}
f(x):=x-\psi(x)  g(\Pi \,I_{\xzero}x)\nu_{\xzero}\,, 
\end{equation*}
where $\nu_{\xzero}:=Q_{\xzero}^{-1}e_n$ is the normal to $J_u$ at $\xzero$, and $\Pi:\R^n\to\R^{n-1}$ is the projection onto the first $n-1$ components.
We compute
\begin{equation*}
 Df=\Id - g \nu_{\xzero}\otimes D\psi - \psi \nu_{\xzero} \otimes (D'\!g\,\Pi\, Q_{\xzero}).
\end{equation*}
Here we use the notation $D'g$ in place of $Dg$ to highlight the derivation in $\R^{n-1}$.
The bounds on $g$ and $\psi$ imply that $|Df-\Id|<3\eps$ everywhere.
In particular, $f$ is a diffeomorphism, which is the identity outside $B_r(\xzero)$.

Let $\mu:=Du\LL B_r(\xzero)- [u](\xzero) \otimes \nu \calH^{n-1} \LL \{I_{\xzero}^{-1}(x', g(x')): x'\in B'_r\}$, where
$\nu$ is the normal to the last set.
By (\ref{eqpropg1}) and (\ref{eqpropg2}), we obtain $|\mu|(\R^n)\le c\eps r^{n-1}$. 

We choose a closed $n-1$-dimensional polyhedron $\hat P$ contained in $B'_{(1-\eps) r}$ and such that $$\calH^{n-1}(B_{(1-\eps) r}'\setminus \hat P)\le \eps
r^{n-1},$$ and define $P_\xzero:=I_{\xzero}^{-1}(\hat P\times \{0\})$ and
\begin{equation}\label{eqdefmup}
 \hat\mu:=D(u\circ f)\LL B_r(\xzero)- [u](\xzero) \otimes \nu_{\xzero} \calH^{n-1} \LL P_\xzero\,.
\end{equation}
By the change-of-variable formula for $BV$ functions, the bounds on $f$ and the estimate in $\mu$ 
we obtain $|\hat\mu|(\R^n)\le c\eps r^{n-1}
\le c\eps|Du|(B_r(\xzero))$. 
All constants may depend only on $n$ and $\calZ$.

 \smallskip
\emph{Step $2$. By a covering argument we conclude the construction.}

Using Vitali's covering theorem, we choose finitely many points $x_1, \dots, x_M\in\R^n$ and radii $r_i\in(0,1)$ with the properties stated in Step 1,
such that 
\begin{equation*}
|Du|\Biggl(\R^n\setminus \bigcup_{i=1}^MB_{r_i}(x_i)\Biggr)<\eps 
\end{equation*}
and the balls $B_{r_i}(x_i)$ are disjoint. Let $f_1,\dots, f_M\in C^1(\R^n;\R^n)$ and 
$P_1,\dots, P_M\subset\R^n$ be the corresponding 
deformations and polyhedra, respectively, and let $u_i^\pm\in \calZ$, $\nu_i\in S^{n-1}$ be the 
corresponding traces and normals.
Let 
$$
f:=f_1\circ f_2\circ \dots\circ f_M\in C^1(\R^n;\R^n).
$$
 Since  $f_i(x)=x$ outside $B_{r_i}(x_i)$ 
we have $$|Df(x)-\Id|+|f(x)-x|\le 6\eps$$ for all $x\in\R^n$. We define $v:=u\circ f$. Then, letting
\begin{equation*}
 \mu^*:=\sum_{j=1}^M (u^+_i-u^-_i)\otimes \nu_i \calH^{n-1}\LL P_i
\end{equation*}
be the polyhedral measure we have constructed in Step 1, we obtain
\begin{align}\nonumber
 |Dv-\mu^* | (\R^n)
 &\le \sum_{i=1}^M |\hat\mu_i|(\R^n)+ |Du|\Biggl(\R^n\setminus \bigcup_{i=1}^M B_{r_i}(x_i)\Biggr)\\
&\le  c\eps|Du|(\R^n)+\eps\,. \label{eqdvmust}
\end{align}
Here $\hat\mu_i$ denotes the analog  for the ball $B_{r_i}(x_i)$
of the remainder $\hat\mu$ obtained in (\ref{eqdefmup}).

\smallskip
\emph{Step $3$. We construct a piecewise-constant $SBV$ function with polyhedral jump set whose gradient is close to $\mu^*$.}

To that end, we will use the results of Lemma \ref{lemmacover} separately stated and proved below.
We consider $c_*$ and the polyhedral decomposition into the cells $\{V_q\}_{q\in G}$ 
of $\R^n$ obtained from Lemma \ref{lemmacover} taking in its hypothesis
the polyhedra $P_1, \dots, P_M$,
with a spacing $\delta>0$ such that $2c_*\delta< \dist(P_i,P_l)$ for all $i\ne l$.

For any $q\in G$, we choose a value $z_q\in \calZ$ such that 
$|V_q\cap v^{-1}(z_q)|=\max_{z'\in \calZ} |V_q\cap v^{-1}(z')|$.
We define $w:\R^n\to\R^N$ by setting
$w(x)=z_q$ if $x\in V_q$.
By the geometric properties of the cells $V_q$ described
in Lemma \ref{lemmacover}, using Poincar\'e's inequality and the trace theorem we
obtain 
\begin{equation}\label{eqpoinccell}
\|v-z_q\|_{L^1(V_q)}\le c \delta |Dv|(V_q) \text{  and }\|v-z_q\|_{L^1(\partial V_q)} \le c |Dv|(V_q)\,,
\end{equation}
where $c$ may depend only on $n$,  $\calZ$ and on $c_*$.
To see this, we observe that  since  $B_{\delta/c_*}(x_q)\subset V_q$ there is $m_q\in\R^N$ such that
$\|v-m_q\|_{L^1(B_{\delta/c_*}(x_q))}\le c \delta |Dv|(V_q)$. 
The estimate is then extended to $V_q$ passing to polar coordinates centered in $x_q$ and using the
one-dimensional Poincar\'e inequality in the radial direction; note that
since $V_q\subset B_{c_*\delta}(x_q)$ the Jacobian determinant is bounded.
Finally one replaces $m_q$ by  $z_q$ using the fact that the volume of the set $V_q\cap\{v=z_q\}$ is at
least $|V_q|/\#\calZ$.
The trace estimate, in turn, follows by using the one-dimensional trace estimate on each segment
connecting a point on $\partial V_q$ with $x_q$, and estimating again the Jacobian determinant using
$B_{\delta/c_*}(x_q)\subset V_q\subset B_{c_*\delta}(x_q)$.

It remains to check that the map $w$ has the desired properties.
Since $w$ takes finitely many values, and is piecewise constant on each of the polyhedra $V_q$ which cover $\R^n$, 
we see that $w\in SBV_\loc(\R^n;\calZ)$ and that $J_w\subset \bigcup_{q\in G}\partial V_q$ is polyhedral.

To estimate $Dw$, we consider two indices
 $q\ne q'\in G$ such that $\calH^{n-1}(\partial V_q\cap \partial V_{q'})>0$.
Denoting by
 $T_qv$ and  $T_{q'}v$  the inner traces of $v$ on the boundaries of $V_q$ and $V_{q'}$ respectively, we obtain,
 using a triangular inequality and (\ref{eqpoinccell}),
\begin{align*}
 |Dw|(\partial V_q\cap \partial V_{q'}) =&
 |z_q-z_{q'}| \calH^{n-1}(\partial V_q\cap \partial V_{q'})\\
\le&
\|T_qv-T_{q'}v\|_{L^1(\partial V_q\cap \partial V_{q'})}\\
&+
\|T_q v-z_{q}\|_{L^1(\partial V_q)}
+\|T_{q'} v-z_{q'}\|_{L^1(\partial V_{q'})}\\
\le& 
|Dv|(\partial V_q\cap \partial V_{q'})
+ c|Dv|(V_q)+c|Dv|(V_{q'})\\
\le& c|Dv|(V_q\cup V_{q'}\cup (\partial V_q\cap \partial V_{q'})).
\end{align*}
If $|\mu^*|(\partial V_q\cap \partial V_{q'})=0$, this estimate suffices. Otherwise,
there is exactly one $j$ such that $\calH^{n-1}(P_i\cap \partial V_q\cap \partial V_{q'})>0$.
Assuming that $\nu_i$ is oriented from $V_{q'}$ to $V_q$, a computation similar to the one above gives
\begin{align*}
 |Dw-\mu^*|(\partial V_q\cap \partial V_{q'}) =&
 \|z_q-z_{q'}-(u_i^+-u_i^-)\chi_{P_i}\|_{L^1(\partial V_q\cap \partial V_{q'})}\\
\le&
\|T_qv-T_{q'}v-(u_i^+-u_i^-)\chi_{P_i}\|_{L^1(\partial V_q\cap \partial V_{q'})}\\
&+
\|T_q v-z_{q}\|_{L^1(\partial V_q)}
+\|T_{q'} v-z_{q'}\|_{L^1(\partial V_{q'})}\\
\le& 
|Dv-\mu^*|(\partial V_q\cap \partial V_{q'})
+ c|Dv|(V_q)+c|Dv|(V_{q'})\\
\le& c|Dv-\mu^*|(V_q\cup V_{q'}\cup (\partial V_q\cap \partial V_{q'})).
\end{align*}
We finally sum over all pairs. 
Since the number of faces of the polyhedra is uniformly bounded, each $V_q$ is 
included only in the estimates for a uniformly bounded number of faces and therefore
\begin{equation*}
 |Dw-\mu^*|(\R^n) \le c  |Dv-\mu^*|(\R^n) \le c\eps |Du|(\R^n)+\eps\,.
\end{equation*}
Since $\eps>0$ was arbitrary, this concludes the proof of Step 3.
\bigskip

The proof of the result then follows by choosing $\eps=1/j$ and defining $f_j$, with a slight abuse of notation,  as the corresponding function $f$ in Step 2. 
\end{proof}

\begin{remark}\rm
  Since in Step 1 of the proof of Theorem \ref{theorem} we may assume $r_y<\eps$, the construction above additionally gives that
  $\dist(x,\supp Du)<{1/j}$ for all $x$ such that $f_j(x)\ne x$.
\end{remark}

\begin{corollary}
In the setting of Theorem \ref{theorem}, if 
 $\psi:S^{n-1}\times \calZ\times\calZ\to[0,\infty)$ is continuous and symmetric and
  $E[u]:=\int_{J_u\cap\Omega} \psi(\nu, u^+, u^-) d\calH^{n-1}$, then $E[u_j]\to E[u]$.
\end{corollary}

\begin{proof}
Since $|D(u\circ f_j)- Du_j|(\Omega)\to0$  we have
$$
\calH^{n-1}(\Omega\cap(J_{u\circ f_j}\triangle J_{u_j}))+
\calH^{n-1}(\{x\in J_{u\circ f_j}\cap J_{u_j}\cap\Omega: u_j^\pm\neq (u\circ f_j)^\pm\})=o(1).
$$
Since $\calZ$ is finite and $S^{n-1}$ compact the function $\psi$ is bounded. Therefore the previous
estimate  implies that 
\begin{equation*}E[u_j]=
\int_{J_{u\circ f_j}\cap \Omega} \psi(\nu_{u\circ f_j}, (u\circ f_j)^+, (u\circ f_j)^-) d\calH^{n-1}+o(1).
\end{equation*}
We remark that $u\circ f_j$ is defined on $f_j^{-1}(\Omega)$, and denote by $\nu_{u\circ f_j}$
the normal to its jump set $J_{u\circ f_j}=f_j^{-1}(J_u)\subset f_j^{-1}(\Omega)$.
One easily checks that $\nu_{u\circ f_j}(x)=Df_j^T(x)\nu(f_j(x))/|Df_j^T(x)\nu(f_j(x))|$.

By the change-of-variables formula, see \cite[Th. 2.91]{AmbrosioFP}, we have
\begin{eqnarray*}
&&\int_{J_{u\circ f_j}\cap \Omega} \psi(\nu_{u\circ f_j}, (u\circ f_j)^+, (u\circ f_j)^-) d\calH^{n-1}
 \\
 &&\qquad\qquad= \int_{J_u\cap f_j(\Omega)} \psi(\nu_j, u^+, u^-) J_{n-1}d^{J_u}f_j^{-1}
 d\calH^{n-1},
\end{eqnarray*}
where $\nu_j:=\nu_{u\circ f_j}\circ f_j^{-1}$ is the normal to $J_{u\circ f_j}$ transported by $f_j$,
which converges uniformly to $\nu$ as $j\to\infty$,
and $ J_{n-1}d^{J_u}f_j^{-1}$ is the Jacobian of the tangential differential of $f_j^{-1}$.
The claim then follows by dominated convergence
using continuity of $\nu\mapsto \psi(\nu, \alpha,\beta)$, that $\nabla f_j$ tends to the identity,
and the fact that $\calH^{n-1}(J_u\setminus f_j(\Omega))\to0$.
\end{proof}

\begin{corollary}\label{Corollario} In the setting of Theorem \ref{theorem}, we obtain that for all $z,z'\in{\cal Z}$ the 
polyhedral sets $A^z_j:=\{x\in\Omega: u_j(x)=z\}$ are such that $\calH^{n-1}(\partial A^z_j\cap \partial A^{z'}_j\cap\Omega)\to
\calH^{n-1}(\partial A^z\cap\partial A^{z'}\cap\Omega)$,
where $A^z:=\{x\in\Omega: u(x)=z\}$ and $\partial$ denotes the reduced boundary.
\end{corollary}
\begin{proof}
It follows from the previous Corollary choosing $\psi(\nu,\alpha,\beta)=1$ if $\{\alpha,\beta\}=\{z,z'\}$ and 
$\psi(\nu,\alpha,\beta)=0$ otherwise.
\end{proof}

We finally state and prove the lemma used in the proof of Step 3 above.

\begin{lemma}\label{lemmacover}
There is $c_*>0$, depending only on $n$, such that the following holds:
 Let $P_1,\dots, P_M$ be $n-1$-dimensional disjoint closed polyhedra in $\R^n$. Then for $\delta>0$
 sufficiently small there are
 countably many pairwise disjoint open convex $n$-dimensional polyhedra $V_q\subset\R^n$, $q\in G$, such that
  $|\R^n\setminus\bigcup_q V_q|=0$
 and $P_j\cap V_q=\emptyset$ for all $j\in\{1, \dots,M\}$ and $q\in G$. 
 For any $q\in G$  there is $x_q\in\R^n$ such that
 $B_{\delta/c_*}(x_q)\subset V_q \subset B_{c_*\delta}(x_q)$. Each polyhedron  $V_q$ has at most $c_*$ faces.
 \end{lemma}
The idea of the proof is to define $G$ as a set of points in $\R^n$ with a spacing of order $\delta$; and then
to construct $(V_q)_{q\in G}$ as the
corresponding Voronoi tessellation. In order for the polyhedral  $P_j$ to be contained in the boundaries
between the 
$V_q$, in a neighbourhood of each $P_j$, we use a grid oriented as $P_j$. The remaining difficulty is to
interpolate between grids of different orientation. This is done superimposing the grids and removing,
in an intermediate layer, some points so that the remaining ones have approximately distance $\delta$ from
each other.
\begin{proof}
We set $\delta_0:=\frac1{5n} \min_{i,j}\dist(P_i, P_j)$. For any $j$ we  define the $t$-neigh\-bour\-hood
of $P_j$ by $(P_j)_t:=\{x\in\R^n: \dist(x,P_j)<t\}$ and fix 
an affine isometry $I_j:\R^n\to\R^n$ such that $P_j\subset I_j(\R^{n-1}\times\{0\})$.

\begin{figure}
 \begin{center}
  \includegraphics[width=10cm]{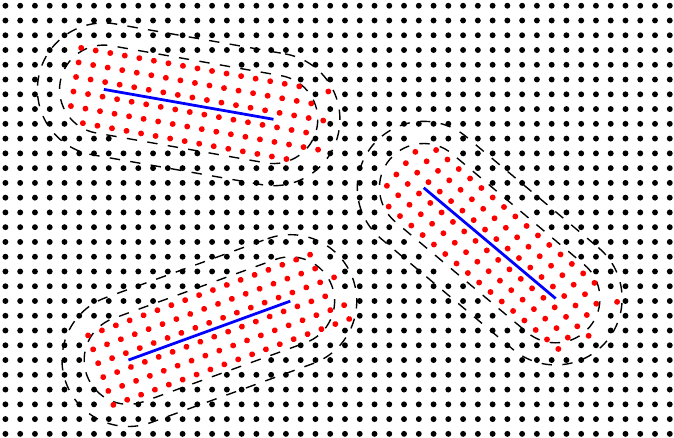}
 \end{center}
\caption{Sketch of the grid construction in the proof of Lemma \ref{lemmacover}.}
\label{fig-grid}
\end{figure}

For $\delta\in (0,\delta_0)$ we set
\begin{equation*}
 \hat G_0:=\delta\Z^n \setminus \bigcup_{j=1}^M (P_j)_{2n\delta_0}
\end{equation*}
and, for $j=1, \dots, M$, 
\begin{equation*}
 \hat G_j:=I_j\Bigl(\delta \Z^n+\frac12\delta e_n\Bigr)\cap (P_j)_{3n\delta_0} \,.
\end{equation*}
The set $\hat G:=\bigcup_{j=0}^M \hat G_j$ is a discrete set with the property that any $x\in\R^n$ has distance
at most  $\sqrt n \delta$ from $\hat G$, see Figure \ref{fig-grid} for an illustration.
Inside each of the disjoint sets $(P_j)_{2n\delta_0}$ the set $\hat G$ coincides with 
$I_j(\delta \Z^n+\frac12\delta e_n)$.
We define $G\subset\hat G$ as a maximal subset with the property that any two points of $G$ have a distance
of at least $\delta/n$. 
By maximality, for any $z\in\hat G\setminus G$ there is $q\in G$ with $|q-z|<\delta/n$; hence
for any $x\in\R^n$ there is a point  $q\in G$ with $|x-q|\le (\sqrt n+1/n)\delta\le n\delta$.
Further,  $G\cap (P_j)_{n\delta_0}=\hat G\cap (P_j)_{n\delta_0}$.

For $q\in G$, let $V_q:=\{x\in\R^n: |x-q|<|x-z|\text{ for all } z\in G, z\ne q\}$. The family of all such  $V_q$ is the Voronoi 
tessellation of $\R^n$ induced by $G$. The $V_q$ are open, disjoint, convex polyhedra which 
cover $\R^n$ up to a null set.
This concludes the construction.

It remains to prove the stated properties.
Since the distance of two points in $G$ is at least $\delta/n$, we have $B_{\delta/(2n)}(q)\subset V_q$.
Since any point in $\R^n$ is at distance smaller than $n\delta$ from a point of $G$, we have
$V_q\subset B_{n\delta}(q)$. In particular, $\overline V_q\cap \overline V_{q'}\ne \emptyset$
implies $|q-q'|\le 2n\delta$. Since the balls $B_{\delta/(2n)}(q)$, with $q\in G$, are disjoint, given $q\in G$ there are at most $(4n^2)^n$ points $q'\in G$ 
such that $|q-q'|\le 2n\delta$.
It follows that $V_q$ is a polyhedron
with at most $(4n^2)^n$ faces.

We finally show that the polyhedra $P_j$ are cointained in the union of the boundaries of the $V_q$.
To do this, fix one $j\in\{1,\dots, M\}$. Set $G_j:=I_j(\delta \Z^n+\frac12\delta e_n)\cap (P_j)_{n\delta_0}$.
By construction, 
$G\cap (P_j)_{n\delta_0}=G_j$. In particular, $P_j\subset \bigcup_{q\in G_j} \overline V_q$.
At the same time, since $G_j$ is symmetric with respect to the hyperplane which contains $P_j$,
each point of $P_j$ is equidistant from at least two of its points, and therefore
 $P_j\subset \bigcup_{q\in G_j} \partial V_q$. This concludes the proof.
\end{proof}

We finally turn to the extension argument which is needed for the derivation of Theorem \ref{theorem} from Theorem \ref{theoremRn}.

  \begin{lemma}\label{lemmaext}
  Let $\Omega\subset\R^n$ be a Lipschitz set with $\partial\Omega$ bounded, let $\calZ\subset\R^N$ be a finite set, and let $u\in SBV_\loc(\Omega;\calZ)$.   Then, there is
  an extension $\tilde u\in SBV_\loc(\R^n;\calZ)$ with  $\tilde u=u$ in $\Omega$,  
 $|D\tilde u|(\partial\Omega)=0$, $|D\tilde u|(\R^n)<c |Du|(\Omega)$.
 \end{lemma}
 \begin{proof}
 To construct the extension, we fix $\eta\in(0,1)$ and $\normal\in C^1(\R^n;\R^n)$ a smoothing of the outer normal  $\nu$
 to $\partial\Omega$, such that $|\normal|=1$ and $\normal\cdot \nu>\eta$ $\calH^{n-1}$-almost everywhere on $\partial\Omega$.
The map $\normal$ is constructed by considering a covering of $\partial\Omega$ by balls in which $\Omega$ is a Lipschitz subgraph, 
in the sense that 
$\Omega\cap B_r(x)=\{y\in B_r(x): (Q_xy)_n<\psi_x(\Pi Q_xy)\}$, with $x\in\partial\Omega$, $Q_x\in O(n)$, $\psi_x:\R^{n-1}\to\R$ Lipschitz,
and $\Pi:\R^n\to\R^{n-1}$ denotes the projection onto the first $n-1$ components.
This implies $Q_x^Te_n\cdot \nu \ge \eta_x:= 1/\sqrt{1+(\mathrm{Lip}(\psi_x))^2}$ on $B_r(x)\cap\partial\Omega$.
By compactness, $\partial\Omega$ is covered by finitely many such balls $\{B_{r_j}(x_j)\}_{j=1,\dots, J}$. We fix 
a partition of unity 
$g_j\in C^\infty_c( B_{r_j}(x_j))$ with $\sum_j g_j=1$ on $\partial\Omega$ and define
$\tildenormal:=\sum_j g_j Q_{x_j}^Te_n$, $\eta:=\min_j \eta_{x_j}$.
It remains only to rescale so that $|\normal|=1$ on $\partial\Omega$. Since 
we already know that $|\tildenormal|\ge \tildenormal\cdot\nu\ge \eta$ on $\partial\Omega$
this can be done setting $\normal:=\varphi(\tildenormal)$, where $\varphi\in C^\infty(\R^n;\R^n)$ coincides with
the projection onto the unit sphere outside $B_\eta(0)$.

Having constructed $\normal$ and $\eta$, we observe that there is $\rho>0$ such that  $(x,t)\mapsto \Phi(x,t):=x + t \normal(x)$ is a
bilipschitz map from $\partial\Omega\times (-\rho,\rho)$ to a tubular neighbourhood of $\partial\Omega$. 
To see this, one first uses  the implicit function theorem on the map $ \R^n\times \R\ni(x,t)\mapsto (\Phi(x,t),t)\in\R^{n}\times\R$ 
to see that it is a diffeomorphism in a neighbourhood of any $(x,t)\in \partial\Omega\times\{0\}$, then the
 compactness of $\partial\Omega$ to show that it is covered by a finite number of such sets, and finally one restricts
 to $x\in\partial\Omega$.
 
 We define $\tilde u(x+t\normal(x))=u(x-t\normal(x))$ for $x\in \partial\Omega$ and $t\in (0, \rho)$,
or equivalently  $\tilde u(x)=u(\Phi(P_t\Phi^{-1}(x))$ for $x\in \Phi(\partial\Omega\times(0,t))$,
where $P_t$ is the linear map that flips the sign of the last argument. We further set
 $\tilde u=u$ in $\Omega$, and $\tilde u$ equal to a constant arbitrary element $z_0$ of $\calZ$
 on the rest of $\R^n$.
 Then $\tilde u:\R^n\to \calZ$.
 By the chain rule for $SBV$ functions, 
 $\tilde u\in SBV_\loc(\R^n;\calZ)$.
By the construction $\tilde u$ has the same trace on both sides of $\partial\Omega$, hence $|D\tilde u|(\partial\Omega)=0$.
  \end{proof}
 \begin{proof}[Proof of Theorem {\rm\ref{theorem}}]
It suffices to apply Theorem \ref{theoremRn} to the extension $\tilde u$ of $u$ constructed in Lemma \ref{lemmaext}. 
 \end{proof}

\begin{remark}\rm
In the statement of Theorem \ref{theoremRn} we can replace the Lipschitz and boundedness assumption on $\Omega$ by the requirement that an extension as in Lemma \ref{lemmaext} 
exists. Such an assumption is satisfied for example if $\Omega$ is a half space, taking the extension by reflection.
\end{remark}

\section*{Acknowledgements}
SC acknowledges financial support by the Deutsche Forschungsgemeinschaft through the Sonderforschungsbereich 1060 {\sl ``The mathematics of emergent effects''}.

\bibliographystyle{alpha-noname}
\bibliography{Bra-Co-Gar}

\end{document}